\documentclass[10pt,reqno]{amsart}
\usepackage[english, activeacute]{babel}
\usepackage{epsfig}
\makeatletter \numberwithin{equation}{section}
\usepackage{amssymb}
\usepackage{latexsym}
\usepackage{amsmath}
\usepackage{bbm}
\usepackage{amsthm}
\usepackage{stmaryrd}
\usepackage{amsfonts}
\date{}
\newcommand{\equ}[1]{(\ref{#1})}
\newcommand{\R}{\mathbb{R}}  

\newcommand{\la}{\lambda}
\newcommand{\dist}{\operatorname{dist}}     
\newcommand{\supp}{\operatorname{supp}}       
\providecommand{\abs}[1]{\lvert#1 \rvert}  
\providecommand{\norm}[1]{\lVert#1 \rVert} 

\newcommand{\ve}{\varepsilon}

\theoremstyle{plain}
\newtheorem{claim}{Claim}

\newtheorem{teo}{Theorem}
\newtheorem{prop}{Proposition}[section]

\newtheorem{lema}{Lemma}[section]
\newtheorem{remark}{Remark}[section]

\theoremstyle{definition}

\begin{document}
\title[ two-dimensional Lazer-McKenna conjecture ]
{ The two-dimensional Lazer-McKenna conjecture  for an exponential
nonlinearity}
\author{Manuel del Pino and Claudio Mu\~noz}
\address{Departamento de Ingenier\'ia Matem\'atica, Universidad de Chile, Casilla 170, Correo 3, Santiago,
Chile.} \email{delpino@dim.uchile.cl, cmunoz@dim.uchile.cl }
\begin{abstract}

We consider the problem of Ambrosetti-Prodi type
\begin{equation}\label{0}\quad\begin{cases}
\Delta u + e^u = s\phi_1 + h(x) &\hbox{ in } \Omega ,\\
 u=0 & \hbox{ on } \partial \Omega,\\
\end{cases}
\nonumber \end{equation} where $\Omega$ is a bounded, smooth domain
in $\R^2$, $\phi_1$ is a positive first eigenfunction of the Laplacian under Dirichlet boundary conditions
and $h\in\mathcal{C}^{0,\alpha}(\bar{\Omega})$. We prove that given $k\ge 1$ this problem has at
least $k$ solutions for all sufficiently large $s>0$, which  answers
affirmatively a conjecture by Lazer and McKenna \cite{LM1} for this
case. The solutions found exhibit multiple concentration behavior
around maxima of $\phi_1$ as $s\to +\infty$.
\end{abstract}

\maketitle

\bigskip
\section{Introduction and statement of main results}

Let $\Omega \subseteq \R^2$ be a bounded and smooth domain. This
paper deals with the  boundary value problem

\begin{equation}\label{In1}  \quad\begin{cases}
\Delta u + e^u = s\phi_1 + h(x) &\hbox{ in } \Omega ,\\
 u=0 & \hbox{ on } \partial \Omega ,\\
\end{cases}
\end{equation}
where $h\in C^{0,\alpha}(\bar{\Omega})$ is given, $s$ is a large,
positive parameter and $\phi _1 $ is a positive first
eigenfunction of of the problem $-\Delta\phi = \la \phi $ under
Dirichlet boundary condition in $\Omega$. We denote its
eigenvalues as
 $$0<\lambda_1<\lambda_2\leq
\lambda_3\leq \cdots  $$

The {\em Ambrosetti-Prodi problem} is the equation

\begin{equation}\label{In1a}  \quad\begin{cases}
\Delta u + g(u) = f(x) &\hbox{ in } \Omega ,\\
 u=0 & \hbox{ on } \partial \Omega ,\\
\end{cases}
\end{equation}
 where $\Omega\subset \R^N$ is bounded and smooth, $f\in \mathcal{C}^{0,\alpha}(\bar{\Omega})$, and the limits
$$
\nu \equiv \displaystyle{\lim_{t\to -\infty} \frac{g(t)}{t}}\ < \
 \mu \equiv \displaystyle{\lim_{t\to+\infty} \frac{g(t)}{t}}
$$
are assumed to exist.  Problem \equ{In1} corresponds to a case in
which $\nu =0$ and $\mu =+\infty$. In 1973, Ambrosetti and Prodi
\cite{AP} assumed that
$$
0< \nu < \la_1 < \mu < \la_2
$$
and additionally that $g''>0$. They proved the existence of a
$C^1$ manifold  ${\mathcal M}$ of codimension 1 which separates
$C^{0,\alpha}(\bar{\Omega})$ into two disjoint open regions,
$$C^{0,\alpha}(\bar{\Omega})= {\mathcal O}_0 \cup {\mathcal M} \cup {\mathcal
O}_2 ,$$ such that Problem  (\ref{In1a}) has no solutions for
$f\in {\mathcal O}_0$, exactly two solutions if $f\in {\mathcal
O}_2$, and exactly one solution if $f\in {\mathcal M}$.

\medskip
In 1975, Berger and Podolak \cite{BePo} obtained a more explicit
representation for the result in \cite{AP} by decomposing
$$f=s\phi_1 + h, \quad \int_{\Omega} h\phi_1 =0 \, ,$$
and proving that for each such an $h$
 there is a  number $\alpha (h)$ such that the  problem

\begin{equation}\label{In1aa} \quad\begin{cases}
\Delta u + g(u) = s\phi_1 + h &\hbox{ in } \Omega ,\\
 u=0 & \hbox{ on } \partial \Omega ,\\
\end{cases}
\end{equation}
has no solution if $s<\alpha(h)$ and  exactly two solutions if
$s>\alpha(h)$. Written in this form, letting $s$ be a parameter
and $h$ fixed, is what is commonly referred to as the
Ambrosetti-Prodi problem.

\medskip
The convexity assumption in the multiplicity result for large and
positive $s$ was relaxed subsequently in \cite{AH,Da1,KW}. In
\cite{LM1}, Lazer and McKenna obtained a third solution of
(\ref{In1aa}) under  the further assumption

\begin{equation*}
\nu <\lambda_1<\lambda_2< \mu <\lambda_3,
\end{equation*}
while  a fourth solution under this circumstance was found by
Hofer \cite{Ho} and by Solimini \cite{So}. In \cite{LM1} it was
further conjectured  that the number of solutions for very large
$s>0$ grows as the interval $(\nu, \mu)$ contains more and more
eigenvalues, in particular, they conjectured that if
\begin{equation}\nu<\la_1 < \mu=+\infty \label{lm}\end{equation}
and $g$ does not grow ``too fast'' at infinity, then for all
$k\ge1$ there is a number $s_k$ such that for all $s>s_k$, Problem
(\ref{In1aa}) has at least $k$ solutions.

Surprisingly enough, Dancer \cite{Da2} was able to disprove the
conjecture in the asymptotically linear case in which $\nu$ and
$\mu$ are finite, exhibiting an example in $N\ge 2$ in which the
interval $(\nu,\mu)$ contains a large number of eigenvalues but
{\em no more than four solutions} for large $s$ exist. The
conjecture, for both $\mu$ finite and infinite actually holds true
in one-dimensional and radial cases under various situations, see
\cite{castrokurepa,hlm,LM2,PMM,RSo} for these and related results.
See also \cite{BMP,DF,DeY,So1}  for other results in the PDE case.

\medskip
How fast should ``too fast" be in the growth of $g$ under the
situation \equ{lm}? The authors of the conjecture had probably in
mind a growth not beyond critical for the nonlinearity. This
constraint was indeed used in \cite{castrokurepa} in the radial
case.

\medskip
Recently Dancer and Yan \cite{DY2,DY3} proved that the
Lazer-McKenna conjecture \emph{holds true}  when $N\ge 3$ and
$$ g(t) = \la t + t_+^p,  \quad 1< p < \frac{N+2}{N-2}, \ \la< \la_1\, , $$
by constructing and describing asymptotic behavior of the
solutions found as $s\to +\infty$. In this case $\nu = \la$ and
$\mu =+\infty$. This has also been done in the critical case
$p={N+2\over N-2}$ if, in addition, $0<\la$ and $N\ge 7$, by Li,
Yan and Yang in \cite{LYY}.

\medskip
Problem (\ref{In1}) is also a problem involving criticality in
$\R^2$. While, strictly speaking, the nonlinearity stays below the
threshold of compactness given by Trudinger-Moser embedding, for
which $e^{u^2}$ is critical, two dimensional equations involving
$e^u$ exhibit {\em bubbling phenomena}, similar to that found at
the critical exponent in higher dimensions. This has been a
subject broadly treated in the literature,
in what regards to construction and classification of unbounded
families of solutions for this type of exponential nonlinearities.

\medskip
The main  result of this paper is  a \emph{positive} answer to the
Lazer-McKenna conjecture for Problem (\ref{In1}).  Given any $m\ge
1$, there are at least $m$ solutions for all $s>0$ sufficiently
large. These solutions can be explicitly described: they exhibit
multiple bubbling behavior around maximum points of $\phi_1$.
\begin{teo}\label{teo1}
Given any $m\ge 1$  and any $s$ sufficiently large, there exists a
solution $u_s$ of Problem $\equ{In1}$ such that
$$
\lim_{s\to+\infty} \int_\Omega e^{u_s} = 8\pi m\, .
$$
More precisely, given any subset $\Lambda$ of $\Omega$ for which
$$ \sup_{\partial\Lambda} \phi_1 <\sup_{\Lambda} \phi_1\, $$ and a
sequence $s\to +\infty$, there is a subsequence and $m$ points
$\xi_i\in \Lambda$ with
$$
\phi_1(\xi_i) = \sup_{\Lambda }  \phi_1
$$
such that as $s\to +\infty$
$$
e^{u_s} \rightharpoonup 8\pi \sum_{i=1}^m \delta_{\xi_i}\, .
$$
\end{teo}

In particular, we observe that  associated to any isolated local
maximum point of $\xi_0$ of $\phi_1$  one has the phenomenon of
multiple bubbling at a single point, namely $ e^{u_s}
\rightharpoonup 8\pi m \delta_{\xi_0}\, $.

\medskip
 The construction gives much more accurate
information on the asymptotic profile of these solutions, in
particular we  have the expansion
$$
u_s = - \frac{s}{\lambda_1}\phi_1 - \rho  + \sum_{i=1}^m G(\xi_i,
x) + o(1)
$$
uniformly on compact subsets of $\bar\Omega\setminus
\{\xi_1,\ldots ,\xi_m\}$, where $\rho= (-\Delta)^{-1}h$ in
$H_0^1(\Omega)$
and $G(x,\xi)$ denotes the symetric Green's function of the problem

\begin{equation}\label{Un1}  \quad\begin{cases}
-\Delta_x G(x,\xi)  =  8\pi\delta_\xi(x), & x \in \Omega ,\\
 \quad G(x,\xi) =  0, & x\in  \partial \Omega. \\
\end{cases}
\end{equation}

In order to restate the problem in perhaps more familiar terms,
let us substitute $u$ in equation \equ{In1} by $u-
\frac{s}{\lambda_1}\phi_1 - \rho$.  Replacing further the
parameter $s$ by $\la_1s$ and setting $k(x) = e^{-\rho}$,
\equ{In1} becomes equivalent to

\begin{equation}\label{In3}  \quad\begin{cases}
\Delta u + k(x)e^{-s\phi_1}e^u=0 &\hbox{ in } \Omega ,\\
 u=0 & \hbox{ on } \partial \Omega,\\
\end{cases}
\end{equation}
and thus what one typically expects  are solutions of $u_s$
\equ{In3} that resemble
$$ u_s (x) \sim \sum_{j=1}^k m_jG(\xi_i, x) , $$
with $m_j >1$, where $\xi_i$'s are maxima of $\phi_1$. This
multiple bubbling phenomenon is in strong opposition to the
seemingly similar, well studied problem
\begin{equation}\label{In3a}  \quad\begin{cases}
\Delta u + \ve^2 k(x)e^u=0 &\hbox{ in } \Omega ,\\
 u=0 & \hbox{ on } \partial \Omega,\\
\end{cases}
\end{equation}
with $k\in\mathcal{C}^2(\bar{\Omega})$,  $\inf_{\Omega}k>0$ and
$\ve\to 0$, where bubbling of solutions with
$$ \int_\Omega  \ve^2 k(x)e^u = O(1) $$ is forced to be  simple, namely with all $m_j$'s equal to one, as it follows from the
results in \cite{BM,LiSha,mawei,NS}. Blowing up families of
solutions to this problem have been constructed  in
\cite{BaPa,ChenLin,KMP,EGP}. For instance it is found in
\cite{KMP} the presence of solutions with arbitrary number of
bubbling points whenever $\Omega$ is not simply connected, see
also \cite{EGM} for a similar phenomenon for large exponents in a
power nonlinearity. Multiple bubbling has been built recently, in
\cite{Wei}, for the anisotropic problem
\begin{equation}  \quad\begin{cases}
{\rm div}\,(a(x)\nabla u )+ \ve^2 k(x)e^u=0 &\hbox{ in } \Omega ,\\
 u=0 & \hbox{ on } \partial \Omega,\\
\end{cases}
\nonumber\end{equation} around isolated local maxima of the
(uniformly positive) coefficient $a$. The moral of our  result is
that multiple bubbling in the isotropic case may be triggered by
the fact that the coefficient in front of $e^u$ does not go to
zero in uniform way. Multiple bubbling ``wants to take place''
where the coefficient vanishes faster in $s$. This should be
somehow connected with phenomena associated to \equ{In3a} where
$k(x)$ is replaced by $|x|^\alpha k(x)$, weight resulting for
Liouville type equations with singular sources. Important advances
in understanding of blowing-up solutions for that problem have
been obtained, see for instance \cite{tarantello} and references
therein.

\medskip
The rest of this paper will devoted to the Proof of Theorem
\ref{teo1}. We will actually give to it a precise version in terms
of Problem \equ{In3} in Theorem \ref{T1} below.

\medskip
As we have mentioned, we do not intend to express our results
in their most general  forms. For instance the choice of $\phi_1$ as the
positive function in the right hand side of \equ{In1} is made for
historical reasons but it is certainly not essential. We could in
principle replace it for instance by any positive function $\phi$,
where now concentration will take place around local maxima of the
function $(-\Delta)^{-1}\phi$ in $H_0^1(\Omega)$.

On the other hand we also remark
that a similar result to Theorem \ref{teo1} is valid for the
problem
\begin{equation}\quad\begin{cases}
\Delta u + \la u+ e^u = s\phi_1 + h(x) &\hbox{ in } \Omega ,\nonumber\\
 u=0 & \hbox{ on } \partial \Omega,\\
\end{cases}
\nonumber \end{equation} provided that $\la <\la_1$. Note that
$\nu =\la$, $\mu =+\infty$ in this case.
The basic fact is that $\Delta +\la$ satisfies
maximum principle. Green's function should consistently be
replaced by the one associated to this operator.

\bigskip

\setcounter{section}{0}
\section{Preliminaries and ansatz for the solution}
\bigskip

In what remains of this paper we fix a set $\Lambda$ as in the statement of Theorem 1.
For notational simplicity we assume
\begin{equation*}
\max_{x\in {\bar{\Lambda} }} \, \phi_1(x)=1.
\end{equation*}
What we will do next is to construct a
reasonably good approximation $U$ to a solution of \equ{In3} which
will have as parameters yet to be adjusted, points $\xi_i$ where
the spikes are meant to take place. As we will see, a convenient
set to select $\xi =(\xi_1,\ldots ,\xi_m) $ is

\begin{equation}
\mathcal{O}_{s} \equiv \Big\{ \mathbf{\xi} \in \bar{\Lambda}^m : 1
-\phi_1(\xi_j)\leq \frac{1}{\sqrt{s}},
 \ \forall\,  j=1,\dots ,m, \hbox{ and } \min_{ i\neq j} \abs{\xi_i - \xi_j}
\geq \frac {1}{s^{\beta}}\Big\}, \label{Os}\end{equation}
where the number $\beta >1$ will be specified later.
We thus fix $\xi \in\mathcal{O}_s$.

\medskip
For numbers $\mu_j>0$, $j=1,\dots ,m$, yet to be chosen,   we define
\begin{equation}\label{Aa3}
u_j(x)=u_{j,s}(x)=\log \frac{8\mu_j^2\delta_j^2}{\big(
\mu_j^2\delta_j^2 + \abs{x-\xi_j}^2 \big)^2} + s \phi_1(\xi_j) - \log k(\xi_j),
\end{equation}
so that $u_j$ solves

\begin{equation}\label{Aa4}
\Delta u + k(\xi_j)\delta_j^2 e^u=0 \hbox{ in } \R ^2, \quad
\int_{\R^2} k(\xi_j)\delta_j^2 e^u = 8\pi,
\end{equation}
where, since we are approximating a solution to \equ{In3}, we naturally choose
\begin{equation}
\delta_j=\delta_j(s)\equiv \exp \big\{-\frac s2
\phi_1(\xi_j)\big\}.
\end{equation}\label{deltaj}

Note that $u_j$ is not zero on the boundary of $\Omega $, so that we add to it a harmonic correction
so that boundary condition is satisfied.
Let $H_j(x)$ be the
 solution of
\begin{equation*}\label{Aa5}  \quad\begin{cases}
\Delta H_j=0 &\hbox{ in } \Omega ,\\
 H_j=-u_j &\hbox{ on } \partial \Omega.
\end{cases}
\end{equation*}
We define our first approximation $U(\xi)$
as
\begin{equation}\label{Aa5a}
U(\xi) \equiv \sum_{j=1}^m U_j,\quad  U_j\equiv  u_j + H_j\, .
\end{equation}
As we will see precisely below, $u_j + H_j \sim G(x,\xi_j)$ where
 $G(x,\xi)$ is the  Green function  defined in \equ{Un1}. Let us consider
 $H(x,\xi)$,  its \emph{regular part}, namely
the solution of
\begin{equation}\label{Un2}  \quad\begin{cases}
-\Delta_x H(x,\xi) = 0 & x\in \Omega ,\\
 H(x,y)=\Gamma(x-y) = -4\log \frac{1}{\abs{x-y}}, & x\in\partial \Omega ,\\
\end{cases}
\end{equation}
so that
 $$G(x,y)=H(x,y)-\Gamma(x-y). $$
While $u_j$ is a good approximation to a solution of \equ{In3}
near $\xi_j$, it is not so much the case for $U$, namely
$$U\,= \, u_j\, + \, (H_j + \sum_{k\ne j} u_k),$$
unless the remainder vanishes at main order near $\xi_j$.  This is
 achieved through  the following precise choice of the
parameters $\mu_k$:
\begin{equation}
\log 8\mu_k^2 = \log k(\xi_j) + H(\xi_k,\xi_k) + \sum_{i\neq
k}G(\xi_i,\xi_k).
\label{muk}\end{equation}
Let us observe in particular that
since $\xi\in \mathcal{O}_{s}$,
\begin{equation}\label{mu2}
\frac1C\leq \mu_k \leq Cs^{2\beta}, \quad \hbox{ for all } k=1,\dots ,m.
\end{equation}
some $C>0$.

\medskip
The following lemma expands $U_j$ in $\Omega$.

\begin{lema}\label{AaLe1} Assume $\xi\in\mathcal{O}_{s}$. Then we have

\begin{equation}\label{Aa6}
H_j(x)=H(x,\xi_j)-\log 8\mu_j^2 + \log k(\xi_j) + O(\mu_j^2\delta_j^2),
\end{equation}
uniformly in $\Omega$, and

\begin{equation}\label{Aa7}
u_j(x)=\log 8\mu_j^2 -\log k(\xi_j) -\Gamma(x-\xi_j) + O(\mu_j^2
s^{2\beta}\delta_j^2),
\end{equation}
uniformly in the region $\abs{x-\xi_j}\geq \frac{1}{2s^{\beta}}$,
so that there,

\begin{equation}\label{Aa8}
U_j(x) = G(x,\xi_j)+O(\mu_j^2 s^{2\beta}\delta_j^2)\, .
\end{equation}
\end{lema}
\begin{proof}
Let us prove (\ref{Aa6}). Define $z(x)=H_j(x)+\log 8\mu_j^2 -\log
k(\xi_j)- H(x,\xi_j) $. Since  $z$ is harmonic we have
\begin{eqnarray*}\label{Aa10}
\max _{\overline{\Omega}} \; \abs{z}  & = &  \max_{\partial \Omega}\; \abs{-u_j + \log 8\mu_j^2 -\log k(\xi_j)-\Gamma(\cdot- \xi_j)} \\
& = & \max_{x\in \partial \Omega} \Big| \log \frac{1}{\abs{x-\xi_j}^4} - \log \frac{1}{\big( \mu_j^2\delta_j^2 + \abs{x-\xi_j}^2 \big)^2} \Big| \\
& = &  O(\mu_j^2\delta_j^2),
\end{eqnarray*}
uniformly in $\Omega$, as $s\to \infty$. Expansion (\ref{Aa7}) is
directly obtained by definition of $u_j$ and $\mu_j$.
\end{proof}

\bigskip
Now, let us  write
\begin{equation}
\delta = \delta(s)=e^{-s/2}, \quad\Omega_s = \delta^{-1}\Omega,
\quad \xi_j=\delta\xi_j'. \label{deltagamma}\end{equation} Then
$u$ solves (\ref{In3}) if and only if $v(y)\equiv u(\delta y)-2s$
satisfies
\begin{equation}\label{Aa11}  \quad\begin{cases}
\Delta v +  q(y,s)e^v =0,  \quad \hbox{ in }\Omega_s ,\\
 v(y) = -2s, \quad y\in  \partial \Omega_s ,\\
\end{cases}
\end{equation}
where $$q(y,s)\equiv k(\delta y)\exp \big\{-s(\phi_1(\delta y) -
1) \big\}.$$

\medskip
Let us define $V(y)=U(\delta y)- 2s$, with $U$ our approximate
solution \equ{Aa5a}. We want to measure the size of the error of
approximation  \begin{equation} R\equiv \Delta V + q(y,s)e^V
.\label{error}\end{equation} It is convenient to do  so in terms
of the following norm.

\begin{equation}\label{Aa12}
\norm{v}_{*}=\sup_{y\in \Omega_s} \Big| \Big[\sum_{j=1}^m
\frac{\gamma_j}{(\gamma_j^2 + \abs{y-\xi_j'}^2)^{3/2}} +\delta^2
\Big]^{-1}v(y)\Big|
\end{equation}
where
\begin{equation}\label{gammaj}
\displaystyle{\gamma_j=\mu_j\delta_j\delta^{-1}=\mu_j\exp\big\{
\frac s2(1-\phi_1(\xi_j))\big\}}.
\end{equation}
\emph{Important facts} in the analysis below are the estimates
\begin{equation}\label{gamma}
\frac 1C \leq \gamma_j \leq Cs^{2\beta}e^{\sqrt{s}/2}\, ,\quad
(\delta\, \gamma_j )\leq Cs^{2\beta}e^{- s/4}\, .
\end{equation}
Here and in what follows, $C$  denotes a generic constant
independent of $s$ or $\xi\in\mathcal{O}_s$.
\begin{lema}\label{AsmLe1}  The error $R$ in $\equ{error}$ satisfies
$$\norm{R}_{*}\,\le \, C\,s^{2\beta+1} e^{-s/4}
\quad\hbox{ as $s\to\infty$.}$$
\end{lema}
\begin{proof} We assume first $\abs{y-\xi_k'}\leq\frac {1}{2s^{\beta}\delta} $, for
some index $k$. We have

\begin{eqnarray*}\label{Asm1}
\Delta V(y)  & = &  -\delta^2\sum_{j=1}^m k(\xi_j)e^{-s\phi_1(\xi_j)}
e^{u_j(\delta y)}= -\sum_{j=1}^m\frac{8\gamma_j^2}{(\gamma_j^2
+\abs{y-\xi_j'}^2 )^2}\\
& = & -\frac{8\gamma_k^2}{(\gamma_k^2
+\abs{y-\xi_k'}^2 )^2} + \sum_{j\neq k} O(\mu_j^2
s^{4\beta}\delta^2\delta_j^2).
\end{eqnarray*}
Let us estimate $q(y,s)e^V(y)$. By (\ref{Aa6}) and the definition
of $\mu_j's$,

\begin{eqnarray*}
H_k(x) & = & H(\xi_k, \xi_k) - \log 8\mu_k^2  + \log k(\xi_j) + O(\mu_k^2 \delta_k^2) + O(\abs{x-\xi_k}) \\
& = &  -\sum_{j\neq k} G(\xi_j,\xi_k)+ O(\mu_k^2\delta_k^2) + O(\abs{x-\xi_k}),
\end{eqnarray*}\label{Asm2}
and if $j\neq k$, by (\ref{Aa8})
\begin{equation*}\label{Asm3a}
U_j(x) = u_j(x) + H_j(x) = G(\xi_j,\xi_k) + O(\abs{x-\xi_k}) +
O(\mu_j^2 s^{2\beta}\delta_j^2).
\end{equation*}
Then

\begin{equation}\label{Asm3b}
H_k(x)+\sum_{j\neq k}U_j(x)=  \sum_j O(\mu_j^2 s^{2\beta}\delta_j^2)
+ O(\abs{x-\xi_k}).
\end{equation}
Therefore,
\begin{eqnarray*}\label{Asm3}
q(y,s)e^{V(y)} & = & q(y,s)\delta^4 \exp \Big\{u_k(\delta y) + H_k(\delta y) + \sum_{j\neq k} U_j(\delta y)\Big\}\\
& = & \frac{8\mu_k^2q(y,s)}{( \gamma_k^2  + \abs{y-\xi_k'}^2 )^2k(\xi_k)} \Big\{ 1+ \sum_{j\neq k} O(\mu_j^2 s^{2\beta}\delta_j^2) + O(\delta \abs{y-\xi_k'}) \Big\} \\
& = & \frac{8\gamma_k^2}{(\gamma_k^2 + \abs{y-\xi_k'}^2 )^2}
\;\Big\{ 1 + O(s\delta\abs{y-\xi_k'}) \Big\}
\end{eqnarray*}
We can conclude that in this region
\begin{eqnarray}\label{Asm4}
|R(y)| \leq  C(m,\Omega) \frac{s\gamma_k^2\delta\abs{y-\xi_k'}}{
(\gamma_k^2 + \abs{y-\xi_k'}^2 )^2} + \sum_{j\neq k}
O(\mu_j^2s^{4\beta}\delta^2\delta_j^2). \nonumber\end{eqnarray} If
$\abs{y-\xi_j'}>\frac {1}{2s^{\beta}\delta} $ for all $j$, using
(\ref{Aa6}), (\ref{Aa7}) and (\ref{Aa8}) we obtain  $$\Delta V =
\sum_j O(\mu_j^2 s^{4\beta}\delta^2\delta_j^2) $$ and
$$q(y,s)e^{V(y)} = O(\delta^4\exp \{ -\sum_{j=1}^m
\Gamma(x-\xi_j)\})=O(s^{4m(m-1)\beta}\delta^4). $$
 Hence,
\begin{eqnarray}\label{Asm5}
R(y) = \sum_jO(s^K\delta^2\delta_j^2) \nonumber\end{eqnarray} for
some $K>0$ so that finally

\begin{eqnarray}\label{Asm6}
\norm{R}_{*} = \sum_{k}O(s\gamma_k\delta) \nonumber\end{eqnarray}
and by estimate (\ref{gamma}) the proof is concluded. \end{proof}

\medskip
Next consider the energy functional associated with (\ref{In3})

\begin{equation}\label{In4}
J_s[u]=\frac 12 \int_{\Omega} \abs{\nabla u}^2 - \int_{\Omega}
k(x)e^{-s\phi_1}e^u.
\end{equation}
We will give an asymptotic estimate of $J_s[U]$, where $U(\xi)$ is
the approximation \equ{Aa5a}. The choice of parameters $\mu_j$ as
in \equ{muk} and computations essentially contained in \cite{KMP}
show that the following expansion holds:

\begin{lema} With the election of $\mu_j $'s given by (\ref{muk}),
\begin{equation}\label{AesTeo1}
J_s[U]=16\pi \sum_{ i\neq j} \log \abs{\xi_i-\xi_j} + 8\pi s
\sum_{j=1}^m \phi_1(\xi_j) + O(1),
\end{equation}
where $O(1)$ is uniform in $\xi\in\mathcal{O}_{s}$.
\end{lema}

\medskip

In the subsequent analysis we will stay in the expanded variable
$y\in\Omega_s$  so that  we will look for solutions of problem
(\ref{Aa11}) in the form $v=V +\psi $, where $\psi $ will
represent a lower order correction. In terms of $\psi $, problem
(\ref{Aa11}) now reads

\begin{equation}\label{Asm7} \quad\begin{cases}
\mathcal{L}(\psi) \equiv \Delta \psi + W \psi = -[R + N(\psi)] & \hbox{ in } \Omega_s, \\
\psi = 0 & \hbox{ on } \partial \Omega_s,
\end{cases}
\end{equation}
where $$ \hbox{ $N(\psi) = W[e^\psi - 1 - \psi]$ and
$W=q(y,s)e^V$, }$$

\medskip
Note that $$W(y)=\sum_{j=1}^m \frac{8\gamma_j^2}{(\gamma_j^2 +
\abs{y-\xi_j'}^2)^2}(1+O(s\delta\abs{y-\xi_j'})) \quad \hbox{for }
y\in\Omega_s ,$$ which can be written in the following way

\begin{lema}\label{AsmLe2}
For $y\in \Omega_s$ and $\xi\in \mathcal{O}_{s}$, $W(y)=O
\big(\delta^2\sum_{j=1}^m \delta_j^2 e^{u_j(\delta y)} \big)$, and
then $\norm{W}_{*}=O(1)$.
\end{lema}

\bigskip


\section{The linearized problem}\label{Lp}

In this section we  develop a solvability theory for the linear
operator $\mathcal{L}$ defined in \equ{Asm7}  under suitable
orthogonality constrains. We consider

\begin{equation}\label{Lp1}
\mathcal{L}(\psi)\equiv \Delta \psi + W(y)\psi,
\end{equation}
where $W(y)$ was introduced in (\ref{Asm7}).  By Lemma
\ref{AsmLe2} the operator $\mathcal{L}$ resembles

\begin{equation}\label{Lp2}
\mathcal{L}_0(\psi)\equiv \Delta \psi + \Big( \delta^2\sum_{j=1}^m
\delta_j^2 e^{u_j} \Big)\psi,
\end{equation}
which is a essentially a superposition of linear operators which, after
translations and dilations,  approach as $s\to\infty$ the operator in $\R^2$

\begin{equation}\label{Lp3}
\mathcal{L}_{*}(\psi)\equiv \Delta \psi +
\frac{8}{(1+\abs{z}^2)^2}\psi,
\end{equation}
namely, equation $\Delta v + e^v =0$ linearized around the radial
solution $v(y)=\log \frac{8}{(1+\abs{y}^2)^2}$. The key fact to develop a satisfactory
solvability theory for the operator $\mathcal{L}$ is the non-degeneracy of $v$ up to  the
natural invariances of the equation under translations and dilations.
In fact, if we set
\begin{eqnarray}\label{Lp4}
Z_{0}(z) & = & \frac{\abs{z}^2-1}{\abs{z}^2 + 1}, \\
Z_{i}(z) & = & \frac{4z_i}{ 1 +\abs{z}^2}, \quad i=1,2,
\end{eqnarray}
the only bounded solutions of $\mathcal{L}_*(\psi)=0$ in $\R^2$ are
linear combinations of $Z_{i}$, $i=0,1,2$; see \cite{BaPa} for a
proof.

\medskip
 We define for $i=0,1,2$ and $j=1,\dots,m$,
$$ Z_{ij}(y)\equiv \frac{1}{\gamma_j} Z_{i}\left(\frac{y-\xi_j'}{\gamma_j} \right),\; i=0,1,2.$$
Additionally, let us consider $R_0$ a large but fixed number and
$\chi$ a radial and smooth cut-off function with $\chi\equiv1$ in
$B(0,R_0)$ and $\chi\equiv 0$ in $B(0,R_0+1)^c$. Let
$$\chi_{j}(y)=\chi(\gamma_j^{-1}\abs{y-\xi_j'}),\quad
j=1,\ldots,m.$$
 Given $h\in L^{\infty}(\Omega_s)$, we consider the
problem of finding a function $\psi $ such that for certain scalars $c_{ij}$ one has

\begin{eqnarray}\label{Lp5}
\begin{cases} & \mathcal{L}(\psi)= h + \sum_{i=1}^2 \sum_{j=1}^m c_{ij}\chi_jZ_{ij}, \quad \hbox{ in } \Omega_s, \\
& \psi=0, \quad \hbox{ on } \partial \Omega_s,\\
& \int_{\Omega_s}\chi_j(y)Z_{ij}\psi =0, \hbox{ for all } i=1,2,\;
j=1,\dots , m.
\end{cases}
\end{eqnarray}


\begin{prop}\label{LpTeo1}
There exist positive constants $s_0>0$ and $C>0$ such that for any $h\in
L^{\infty}(\Omega_s)$ and any $\xi\in \mathcal{O}_{s}$, there is a
unique solution $\psi =T(h)$ to problem (\ref{Lp5}) for all
$s>s_0$, which defines a linear operator of $h$. Besides, we
have the  estimate

\begin{equation}\label{LpTeo1a}
\norm{T(h)}_{\infty} \leq C\,s\,\norm{h}_{*}.
\end{equation}

\end{prop}

The proof will be split into  a series of lemmas which we state and prove next.

\begin{lema}\label{LpLe1}
The operator $\mathcal{L}$ satisfies the {maximum principle} in
$\Omega_R \equiv \Omega_s \backslash \cup_{j=1}^m
B(\xi_j',R\gamma_j) $, for $R$ large but independent of $s$. Namely,
if $\mathcal{L}(\psi)\leq 0$ in $\Omega_R$ and $\psi\geq 0$ on
$\partial \Omega_R$, then $\psi\geq 0$ in $\Omega_R$.
\end{lema}

\begin{proof} Notice that for $s$ sufficiently large, $\gamma_j\leq\delta^{-1}$, for all $j$.
This ensures that $\Omega_R$ is well defined. Now, it is
sufficient to find a smooth function $f(y)$ such that $f>0$ in
$\overline{\Omega}_R$ and $\mathcal{L}(f)\leq 0$ in
$\overline{\Omega}_R$.

For this purpose, we use the following lemma, whose proof is
contained in \cite{Wei}:

\begin{lema}\label{LpLe1a}
There exist constants $R_1>0$, $C>0$ such that for any $s>0$ large
enough, there exists $f:\Omega_{R_1}\to[1,\infty)$ smooth and
positive verifying $$\mathcal{L}(f)\leq-\sum_{j=1}^m
\frac{\gamma_j}{\abs{y-\xi_j'}^3} -\delta^2$$ in $\Omega_{R_1}$,
and $1<f\leq C$ uniformly in $\Omega_{R_1}$.
\end{lema}

We briefly recall the argument: we consider numbers $R_1$, $s$
large enough and define
$$\frac 1C \alpha^2 f(y)=f_0(\delta y)- \sum_{j=1}^m \frac{\gamma_j^{\alpha}}{\abs{y-\xi_j'}^{\alpha}},$$
with $f_0$ the solution of $-\Delta f_0 =1$ in $\Omega$, $f_0=2$
on $\partial \Omega$ and $\alpha\in(0,1)$. It is directly checked
that $f$ verifies the required conditions.\end{proof}

Let us consider now the {\em inner norm}

\begin{equation*}\label{Lp10}
\norm{\psi}_{i} \equiv \sup_{\Omega_R^c} \; \abs{\psi}
\end{equation*}
where we understand  $\Omega^c_R \equiv \Omega_s \backslash \Omega_R
= \cup_{j=1}^m B(\xi_j',R\gamma_j)$.

\begin{lema}\label{LpLe2}
There exists a constant $C=C(R,m)>0$ such that if $\mathcal{L}(\psi)
=h$ in $\Omega_s $, $\psi=0$ on $\partial\Omega_s$, $h\in L^{\infty}(\Omega_s)$, and $s$ is
sufficiently large, we have
\begin{equation}\label{Lp11}
\norm{\psi}_{\infty} \leq C \big\{ \norm{\psi}_{i} + \norm{h}_{*}
\big\}.
\end{equation}
\end{lema}
\begin{proof} We will establish this estimate with the aid of Lemmas~\ref{LpLe1} and~\ref{LpLe1a}. We let $f$
be the function defined in the latter result. We consider  the
function $$\hat{\psi } = (\norm{\psi}_{i} + \norm{h}_{*})f, $$ and
claim that $\hat{\psi}\geq \abs{\psi }$ on $\partial \Omega_R$ if
$R$ is sufficiently large. In fact, if $y\in \partial \Omega_s$,
by the positivity of $f$, we have

\begin{equation*}\label{Lp14}
\hat{\psi}(y) \geq 0 = \abs{\psi (y)}.
\end{equation*}
On the other hand, if $\abs{y-\xi_k'}=R\gamma_k$ for some
$k=1,\dots,m$,

\begin{equation*}\label{Lp16}
\hat{\psi}(y) \geq  \norm{\psi}_i f \geq \norm{\psi}_i \geq
\abs{\psi(y)}, \quad \hbox{ for } \abs{y-\xi_k'}=R\gamma_k, \;
k=1,\dots,m.
\end{equation*}
Finally, using that $\displaystyle{\abs{h(y)}\leq \Big(\sum_{j=1}^m
\frac{\gamma_j}{(\gamma_j^2 + \abs{y-\xi_j'}^2)^{3/2}}+\delta^2
\Big)\norm{h}_{*}}$, we have for $y\in \Omega_R$,

\begin{eqnarray*}\label{Lp17}
\mathcal{L}(\hat{\psi})(y) & \leq & (\norm{\psi}_i + \norm{h}_{*})\mathcal{L}(f) \leq - \norm{h}_{*}\Big\{ \sum_{j=1}^m \frac{\gamma_j}{\abs{y-\xi_j'}^3} + \delta^2 \Big\}\\
& \leq & -\norm{h}_{*} \Big\{\sum_{j=1}^m  \frac{\gamma_j}{(\gamma_j^2 +
\abs{y-\xi_j'}^2)^{3/2}} +\delta^2\Big\} \\
&\leq & -\abs{h(y)} \leq -\abs{\mathcal{L}(\psi)(y)}
\end{eqnarray*}
provided $R$ large. In particular, we have
$\mathcal{L}(\hat{\psi})\leq -\mathcal{L}(\psi)$ and
$\mathcal{L}(\hat{\psi})\leq \mathcal{L}(\psi)$, in $\Omega_R$.
Hence, by Maximum Principle in Lemma~\ref{LpLe1} we have $\abs{\psi
(y)} \leq \hat{\psi}(y)$, for $y\in \Omega_R$.
From this we obtain
\begin{equation*}\label{Lp18}
\norm{\psi}_{\infty} \leq \norm{\hat{\psi }}_{\infty} \leq C
\big\{ \norm{\psi}_{i} + \norm{h}_{*} \big\}
\end{equation*}
as desired. \end{proof}

\medskip
The   next   step  is  to  obtain a priori estimates  for the
problem

\begin{eqnarray}\label{Lp19}
\begin{cases} & \mathcal{L}(\psi)= h, \quad \hbox{ in } \Omega_s, \\
& \psi=0, \quad \hbox{ on } \partial \Omega_s,\\
& \int_{\Omega_s}\chi_jZ_{ij}\psi =0, \hbox{ for all } i=0,1,2,\;
j=1,\dots , m.
\end{cases}
\end{eqnarray}
which involves more orthogonality conditions than those in
\equ{Lp5}.
 We have the following estimate.
\begin{lema}\label{LpLe3}
Let  $\psi$ be a solution of Problem $(\ref{Lp19})$ with $\xi\in
\mathcal{O}_{s}$. Then, there exists a  $C>0$ such that
\begin{equation}
\norm{\psi}_{\infty}\leq C\, \norm{h}_{*}
\end{equation}
for all $s>0$ sufficiently large.
\end{lema}

\begin{proof}
We carry out the proof by a contradiction argument. If the result was false,
then, there would exist a sequence $s_n\to
\infty$, points $\xi^n\in \mathcal{O}_{s_n}$, functions $h_n$ with
$\norm{h_n}_{*}\to 0$ and associated solutions $\psi_n$ with
$\norm{\psi_n}_{\infty}=1$ such that
\begin{eqnarray}\label{Lp20}
\begin{cases} & \mathcal{L}(\psi_n)= h_n, \quad \hbox{ in } \Omega_{s_n}, \\
& \psi_n=0, \quad \hbox{ on } \partial \Omega_{s_n},\\
& \int_{\Omega_{s_n}}\chi_jZ_{ij}\psi_n =0, \hbox{ for all }
i=0,1,2,\; j=1,\dots , m.
\end{cases}
\end{eqnarray}
By virtue of Lemma~\ref{LpLe2} and $\norm{\psi_n}_{\infty}=1$ we
have $\liminf_{n\to \infty} \norm{\psi_n}_{i}\geq\alpha >0$. Let
us set $\hat{\psi}_n(z)=\psi_n((\xi_j')^n + \gamma_j^n z)$, where
the index $j=j(n)$ is such that $\sup_{B(\xi_j'^n,R\gamma_j)}
\abs{\psi_n}\geq \alpha$, and can be assumed to be the same for
all $n$. We notice that $\hat{\psi}_n$ satisfies
\begin{equation*}
\Delta \hat{\psi}_n +  (\gamma_j^n)^2 W \, \hat{\psi}_n =
(\gamma_j^n)^2 h_n, \quad \hbox{ in } \Omega_n \equiv
\gamma_j^{-1}(\Omega_s-(\xi_j')^n)\ .
\end{equation*}
Elliptic estimates allow us to assume   that $\hat{\psi}_n$
converges uniformly over compact subsets of $\R^2$ to a bounded,
non-zero solution $\hat{\psi}$ of
\begin{equation*}
\Delta \psi + \frac{8}{(1+\abs{z}^2)^2}\psi =0.
\end{equation*}
This implies that $\hat{\psi}$ is a linear combination of the
functions $Z_{i},\,i=0,1,2$, namely, $\hat{\psi}=\sum_{k=0}^2
\alpha_k Z_{k}$. But orthogonality conditions over $\hat{\psi}_n$
pass to the limit thanks to $\norm{\hat{\psi}_n}_{\infty}\leq 1$.
Dominated convergence then yields

\begin{eqnarray*}
0=\int_{\Omega_{s_n}}\chi_j Z_{ij}\psi_n & = & \int_{\R^2} \chi Z_{i}\hat{\psi}_n + o(1) \\
& = &\sum_{k=0}^2 \alpha _k\int_{\R^2}\chi Z_{i}Z_{k} + o(1), \quad
i=0,1,2.
\end{eqnarray*}
But $\int_{\R^2}\chi Z_{i}Z_{k}=0$ for $i\neq k$ and
$\int_{\R^2}\chi Z_{i}^2>0$. Then $\alpha_k=0$ for all $k=0,1,2$ and
hence $\hat{\psi}\equiv 0$, a contradiction with $\lim\inf_{n\to
\infty}\norm{\psi_n}_i>0$.
\end{proof}

\medskip
Now we will deal with problem (\ref{Lp19}) lifting  the
orthogonality constraints $\int_{\Omega_s}\chi_jZ_{0j}\psi
=0,\quad j=1,\dots ,m$, namely

\begin{eqnarray}\label{Lp21}
\begin{cases} & \mathcal{L}(\psi)= h, \quad \hbox{ in } \Omega_s, \\
& \psi=0, \quad \hbox{ on } \partial \Omega_s,\\
& \int_{\Omega_s}\chi_jZ_{ij}\psi =0, \hbox{ for all } i=1,2,\;
j=1,\dots , m.
\end{cases}
\end{eqnarray}
We have the following a priori estimates for this problem.


\begin{lema}\label{LpLe4} Let $\psi$ be a solution of (\ref{Lp21})
with $\xi\in \mathcal{O}_{s}$. Then, there exists a  $C>0$ such
that
\begin{equation}\label{LpLe4a}
\norm{\psi}_{\infty}\leq C\, s\, \norm{h}_{*}
\end{equation}
for all $s$ sufficiently large.
\end{lema}

\begin{proof}
Let $R>R_0+1$ be a large and fixed number. Let us consider the function
\begin{equation}
\hat{Z}_{0j}=Z_{0j}(y)-\frac{1}{\gamma_j} + a_{0j}G(\delta y,\xi_j),
\end{equation}
where
\begin{equation}\label{Lp22}
a_{0j}\equiv \frac{1}{\gamma_j\{ H(\xi_j,\xi_j) -4\log (\delta
\gamma_jR)\}}\, .
\end{equation}
From estimate  (\ref{mu2}), we have  $$C_1\abs{\log\delta_j} \leq
\log (\delta\gamma_jR)\leq C_2\abs{\log \delta_j} $$ and

\begin{equation}\label{Lp22a}
\hat{Z}_{0j}(y)=O\left(\frac{G(\delta
y,\xi_j)}{\gamma_j\abs{\log\delta_j}}\right).
\end{equation}
Next we consider radial smooth cut-off functions $\eta_1$ and
$\eta_2$ with the following properties:

\begin{eqnarray*}
& & 0\leq\eta_1\leq 1, \quad \eta_1\equiv 1 \hbox{ in } B(0,R), \quad \eta_1\equiv 0 \hbox{ in } B(0,R+1)^c; \hbox{ and }\\
& & 0\leq\eta_2\leq 1, \quad \eta_2\equiv 1 \hbox{ in } B(0,1),
\quad \eta_2\equiv 0 \hbox{ in } B(0,\frac 43)^c.
\end{eqnarray*}
With no loss of generality we assume that $B(0,\frac
43)\subseteq\Omega$. Then we set

\begin{equation}\label{Lp23}
\eta_{1j}(y)=\eta_1\left(\frac{\abs{y-\xi_j'}}{\gamma_j}\right),\quad
\eta_{2j}(y)=\eta_2\left(4\delta\abs{y-\xi_j'}\right),
\end{equation}
and define the test function

\begin{equation*}\label{Lp24}
\tilde{Z}_{0j}=\eta_{1j}Z_{0j} + (1-\eta_{1j})\eta_{2j}\hat{Z}_{0j}.
\end{equation*}

Let $\psi$ be a solution to problem (\ref{Lp21}). We will modify
$\psi$ so that the extra orthogonality conditions with respect to
$Z_{0j}$'s hold. We set

\begin{equation}\label{Lp240}
\tilde{\psi}=\psi + \sum_{j=1}^m d_j \tilde{Z}_{0j} +
\sum_{i=1}^2\sum_{j=1}^m e_{ij}\chi_j Z_{ij}.
\end{equation}
We adjust $\tilde{\psi}$ to satisfy the orthogonality condition

\begin{equation}\label{Lp24a}
\int_{\Omega_s} \chi_jZ_{ij}\tilde{\psi}=0, \quad \hbox{ for all }
i=0,1,2; \; j=1,\dots,m.
\end{equation}
Then,
\begin{equation}\label{Lp25}
\mathcal{L}(\tilde{\psi}) = h + \sum_{j=1}^m
d_j\mathcal{L}(\tilde{Z}_{0j}) + \sum_{i=1}^2\sum_{j=1}^m
e_{ij}\mathcal{L}(\chi_j Z_{ij})\, .
\end{equation}
If (\ref{Lp24a}) holds, the previous lemma allows us to conclude

\begin{equation}\label{Lp26}
\norm{\tilde{\psi}}_{\infty} \leq  C \Big\{ \norm{h}_* +
\sum_{j=1}^m \abs{d_j}\norm{\mathcal{L}(\tilde{Z}_{0j})}_{*} +
\sum_{i=1}^2\sum_{j=1}^m
\abs{e_{ij}}\norm{\mathcal{L}(\chi_j\tilde{Z}_{ij})}_{*}\Big\}.
\end{equation}
Estimate (\ref{LpLe4a}) is a direct consequence of the following two
claims:

\begin{claim}\label{Cl1}The constants $d_j$ and $e_{ij}$ are well defined and
\begin{equation}\label{Lp27}
\norm{\mathcal{L}(\chi_jZ_{ij})}_* \leq \frac{C}{\gamma_j},
\quad\quad \norm{\mathcal{L}(\tilde{Z}_{0j})}_* \leq
\frac{C\log s}{\gamma_j\abs{\log \delta_j}}, \quad i=1,2;\;
j=1,\dots,m.
\end{equation}
\end{claim}
\begin{claim}\label{Cl2} The following bounds hold.
\begin{equation}\label{Lp28}
\abs{d_j}\leq C\gamma_j \abs{\log\delta_j}\, \norm{h}_*, \quad \abs{e_{ij}}\leq C\gamma_j \log s\,
\norm{h}_*, \quad i=1,2;\; j=1,\dots ,m.
\end{equation}
\end{claim}

\medskip
After these facts have been established,  using that
$$\hbox{ $\displaystyle{\norm{\tilde{Z}_{0j}}_{\infty}\leq
\frac{C}{\gamma_j}}$ and $\displaystyle{\norm{\chi_j
Z_{ij}}_{\infty}\leq \frac{C}{\gamma_j}} $, }$$ we obtain
(\ref{LpLe4a}), as desired.

\medskip
 Let us prove now Claim~\ref{Cl1}. First we find $d_j$ and $e_{ij}$. From definition (\ref{Lp240}),
 orthogonality conditions (\ref{Lp24a}) and  the fact that $\supp
\chi_j \chi_k =\emptyset$ if $j\neq k$, we can write

\begin{equation}\label{App0}
e_{ij}=-\frac{\sum_{k=1}^m d_k \int_{\Omega_s} \chi_jZ_{ij}\tilde{Z}_{0k}}{\int_{\Omega_s} \chi_jZ_{ij}^2}, \quad i=1,2;\, j=1,\dots,m.
\end{equation}
Notice that $\displaystyle{\int_{\Omega_s} Z^2_{ij}\chi_j^2}
=c>0$, for all $i,j$, and $$\int_{\Omega_s}
\chi_jZ_{ij}\tilde{Z}_{0l}=O\left( \frac{\gamma_j \log s}{\gamma_l
\abs{\log\delta_l}}\right),\quad j\neq l.$$ Then, from
(\ref{App0})
\begin{equation}\label{App1}
\abs{e_{ij}}\leq C\sum_{l\neq j} \abs{d_l}\frac{\gamma_j \log s}{\gamma_l \abs{\log \delta_l}}.
\end{equation}
We need to show that $d_j$ is well defined. In fact, multiplying definition (\ref{Lp240}) by $Z_{0k}\chi_k$, integrating and using the orthogonality condition (\ref{Lp24a}) for $i=0$, we get
\begin{equation}\label{App2}
\sum_{j=1}^m d_j \int_{\Omega_s}\chi_k Z_{0k} \tilde{Z}_{0j} = - \int_{\Omega_s} \chi_k Z_{0k}\psi,\quad \forall k=1,\dots,m.
\end{equation}
But $\displaystyle{\int_{\Omega_s} \chi_k Z_{0k} \tilde{Z}_{0k}=\int_{\Omega_s}\chi_k Z_{0k}^2 =C}$, for all $k$, and
 $$\displaystyle{\int_{\Omega_s} \chi_k Z_{0k} \tilde{Z}_{0j}=O\left( \frac{\gamma_k \log s}{\gamma_j \abs{\log \delta_j}}\right)},$$
if $k\neq j$. Then, if we define
\begin{equation*}
m_{kj}=\int_{\Omega_s}\chi_k Z_{0k} \tilde{Z}_{0j}, \hbox{ and } f_k=- \int_{\Omega_s} \chi_k Z_{0k}\psi,
\end{equation*}
system (\ref{App2}) can be written as $$\sum_{j=1}^m  m_{kj} d_j=
f_k, \qquad k=1,\dots,m$$ But the matrix with coefficients
$\gamma_j m_{kj} \gamma_k^{-1}$ is clearly  diagonal-dominant,
thus invertible, so  the  matrix $m_{kj}$ is also invertible. Thus
$d_k$ is well defined.

\medskip
Let us prove  inequalities \equ{Lp27}. We note that in the region
$\abs{y-\xi_j'}\leq (R+1)\gamma_j$,
\begin{eqnarray*}
\mathcal{L}(\chi_jZ_{ij})& = &\chi_j\mathcal{L}(Z_{ij}) + \Delta\chi_jZ_{ij} + 2\nabla\chi_j\cdot \nabla Z_{ij}\\
& = & O\left( \frac{s\delta_j\gamma_j^2}{(\gamma_j^2 + \abs{y-\xi_j'}^2)^2} \right) + O\left( \frac{\gamma_j^{-1}}{(\gamma_j^2 + \abs{y-\xi_j'}^2)^{1/2}} \right)\\
& & \quad\quad + O\left( \frac{1}{(\gamma_j^2 + \abs{y-\xi_j'}^2)^{3/2}} \right),
\end{eqnarray*}
and then $\norm{\mathcal{L}(\chi_jZ_{ij})}_{*}=O(\gamma_j^{-1})$. We prove now the second inequality in (\ref{Lp27}). In fact,

\begin{equation*}
\mathcal{L}(\tilde{Z}_{0j}) = \Delta \eta_{1j}(Z_{0j}-\hat{Z}_{0j})+
2\nabla\eta_{1j}\cdot \nabla(Z_{0j}-\hat{Z}_{0j})+
2\nabla\eta_{2j}\cdot\nabla \hat{Z}_{0j}+
\Delta\eta_{2j}\hat{Z}_{0j} \, +
\end{equation*}
\begin{equation*}
+\,\,\eta_{1j}\big\{\mathcal{L}(Z_{0j})-\mathcal{L}(\hat{Z}_{0j})
\big\} + \eta_{2j}\mathcal{L}(\hat{Z}_{0j}),
\end{equation*}
Now we consider the four regions $$\Omega_{1}\equiv \left\{
\abs{y-\xi_j'}\leq \gamma_jR\right\}, \quad
\Omega_{2}\equiv \left\{ \gamma_jR <
\abs{y-\xi_j'}\leq \gamma_j(R+1)\right\},$$
$$\Omega_{3}\equiv \left\{ \gamma_j (R+1) <
\abs{y-\xi_j'}\leq \frac{1}{4\delta} \right\}, \hbox{
and } \Omega_{4}\equiv \left\{ \frac{1}{4\delta} <
\abs{y-\xi_j'}\leq \frac{1}{3\delta}\right\}.$$
Notice that (\ref{Asm4}) and (\ref{Lp3}) imply
\begin{equation}\label{1}
\eta_{1j}\big\{\mathcal{L}(Z_{0j})-\mathcal{L}(\hat{Z}_{0j}) \big\}
+ \eta_{2j}\mathcal{L}(\hat{Z}_{0j})=O\left(
\frac{s\delta\gamma_j^2}{(\gamma_j^2 + \abs{y-\xi_j'}^2)^{3/2}}
\right)
\end{equation}
for all $y\in \Omega_{1}\cup\Omega_{2}$.
Lut us now analyze $\mathcal{L}(\tilde{Z}_{0j})$ in each $\Omega_{i}$. In $\Omega_{1}$,

\begin{equation}\label{2}
\mathcal{L}(\tilde{Z}_{0j})=O\left(
\frac{s\delta\gamma_j^2}{(\gamma_j^2 +
\abs{y-\xi_j'}^2)^{3/2}}\right).
\end{equation}
In $\Omega_{2}$,

\begin{equation}\label{App3}
Z_{0j}-\hat{Z}_{0j}=\frac{1}{\gamma_j}-a_{0j}G(\delta y ,\chi_j)=-a_{0j}\left\{ 4\log\frac{\gamma_jR}{\abs{y-\xi_j'}}+ O(\delta\gamma_j)\right\},
\end{equation}
hence we conclude

\begin{equation}\label{App4}
\abs{Z_{0j}-\hat{Z}_{0j}}=O\left(
\frac{1}{\gamma_j\abs{\log\delta_j}}\right),\quad \hbox{ and } \left|\nabla
(Z_{0j}-\hat{Z}_{0j})\right|=O\left(
\frac{1}{\gamma_j^2\abs{\log\delta_j}}\right),
\end{equation}
and then
\begin{equation}\label{3}
\mathcal{L}(\tilde{Z}_{0j})=O\left(
\frac{1}{\gamma_j^2\abs{\log\delta_j}}\right).
\end{equation}
In $\Omega_{4}$, thanks to (\ref{Lp22a}), $\displaystyle{
\abs{\hat{Z}_{0j}}=O\left(\frac{1}{\gamma_j\abs{\log \delta_j}}\right)}$,
$\displaystyle{\abs{\nabla\hat{Z}_{0j}}=O\left(
\frac{\delta}{\gamma_j\abs{\log\delta_j}}\right)}$ and
\begin{eqnarray*}
\mathcal{L}(\hat{Z}_{0j})& = & \Delta Z_{0j} + W\hat{Z}_{0j}\\
& = & O\left(\frac{\gamma_j}{(\gamma_j^2 + \abs{y-\xi_j'}^2)^2} \right)+ \sum_{k\neq j}O\left( \frac{1}{\gamma_j\abs{\log\delta_j}} \frac{\gamma_k^2}{(\gamma_k^2 + \abs{y-\xi_k'}^2)^2} \right) \, .\\
\end{eqnarray*}
Then, in this region
\begin{equation}\label{4}
\norm{\mathcal{L}(\tilde{Z}_{0j})}_*=O\left( \frac{1}{\gamma_j\abs{\log\delta_j}}\right).
\end{equation}
Finally, we consider $y\in\Omega_{3}$. We have

\begin{eqnarray*}
\mathcal{L}(\tilde{Z}_{0j})& = & \mathcal{L}(\hat{Z}_{0j})\\
& = & \left\{ W-\frac{8\gamma_j^2}{(\gamma_j^2 + \abs{y-\xi_j'}^2)^2}\right\}Z_{0j} + W\left\{ a_{0j}G(\delta y ,\xi_j) -\frac{1}{\gamma_j}\right\} \\
& \equiv  & A_1 + A_2.
\end{eqnarray*}
To estimate these two terms, we need to split $\Omega_{3}$ into
several subregions. We let $$\Omega_{3,j}\equiv
\displaystyle{\left\{ \gamma_j(R+1)< \abs{y-\xi_j'}\leq
\frac{1}{2s^{\beta} \delta} \right\}},$$ $$ \Omega_{3,k}\equiv
\left\{ y\in\Omega_3 \;\Big|\; \abs{y-\xi_k'}\leq
\frac{1}{2s^{\beta} \delta} \right\},\;k\neq j,$$
$$ \hbox{ and } \tilde{\Omega}_{3}\equiv \left\{ y\in\Omega_3 \;\Big|\;
\abs{y-\xi_l'}\geq \frac{1}{2s^{\beta} \delta},\ \forall \,l
\right\}.$$  From Lemma \ref{AsmLe1},
$\displaystyle{A_1=O\left(\frac{s\gamma_j\delta}{(\gamma_j^2 +
\abs{y-\xi_j'}^2)^{3/2}}\right)}$ in $\Omega_{3,j}$, and
$A_1=O(s^{K}\delta^2\delta_j^2\gamma_j^{-1})$ in
$\tilde{\Omega}_3$.

\medskip
If $y\in\Omega_{3,j}$,
\begin{eqnarray*}
A_2 & = & O\left(\frac{\gamma_j^2 a_{0j}}{(\gamma_j^2 + \abs{y-\xi_j'}^2)^2}\Big\{-\log\gamma_jR +\log\abs{y-\xi_j'}+ \delta\abs{y-\xi_j'}\Big\}\right)\\
& = & O\left(\frac{1}{\abs{\log\delta_j}}\frac{1}{(\gamma_j^2 + \abs{y-\xi_j'}^2)^{3/2}} \right),
\end{eqnarray*}
and $A_2=O(s^K\delta^2\delta_j^2)$, for some large $K$. Finally we get, for all $y\in\Omega_{3,j}\cup\tilde{\Omega}_3$,

\begin{equation}\label{5}
\abs{\mathcal{L}(\tilde{Z}_{0j})}=O\left( \frac{1}{\abs{\log\delta_j}}\frac{1}{(\gamma_j^2 + \abs{y-\xi_j'}^2)^{3/2}} \right).
\end{equation}

In $\Omega_{3,k}$, $k\neq j$, we write
\begin{eqnarray*}
\mathcal{L}(\tilde{Z}_{0j}) & = & \Delta Z_{0j} + W \hat{Z}_{0j}\\
& = & \frac{-8\gamma_j^2}{(\gamma_j^2 + \abs{y-\xi_j'}^2)^2}Z_{0j} + W\hat{Z}_{0j}\\
& = & O( s^K \delta_j\delta^3) + O\left(\frac{\gamma_k^2}{(\gamma_k^2 + \abs{y-\xi_k'}^2)^{2}}\frac{G(\delta y,\xi_j)}{\gamma_j\abs{\log\delta_j}} \right)\\
& = & O\left( \frac{\gamma_k^2}{(\gamma_k^2 + \abs{y-\xi_k'}^2)^{2}}\frac{\log s}{\gamma_j\abs{\log\delta_j}}\right),
\end{eqnarray*}
and then, combining (\ref{1})-(\ref{5}) and the previous estimate, we arrive at

$$\norm{\mathcal{L}(\tilde{Z}_{0j})}_{*}=O\left(
\frac{\log s}{\gamma_j\abs{\log \delta_j}}\right).$$ Finally, we
prove Claim \ref{Cl2}. Testing equation (\ref{Lp25}) against
$\tilde{Z}_{0j}$ and using relations (\ref{Lp26}), (\ref{Lp27}),
we get

\begin{eqnarray*}
& & \sum_{k=1}^m d_k\int_{\Omega_s}\mathcal{L}(\tilde{Z}_{0k})\tilde{Z}_{0j} =\\
& = & -\int_{\Omega_s}h\tilde{Z}_{0j}-\int_{\Omega_s}\tilde{\psi}\mathcal{L}(\tilde{Z}_{0j}) + \sum_{l=1}^2\sum_{k=1}^m e_{lk}\int_{\Omega_s}\chi_kZ_{lk}\mathcal{L}(\tilde{Z}_{0j})\\
& \leq & C\frac{\norm{h}_*}{\gamma_j} + C\norm{\tilde{\psi}}_{\infty}\norm{\mathcal{L}(\tilde{Z}_{0j})}_* + C\sum_{l=1}^2\sum_{k=1}^m \abs{e_{lk}}\frac{\norm{\mathcal{L}(\tilde{Z}_{0j})}_{*}}{\gamma_k}
\end{eqnarray*}
\begin{eqnarray*}
& \leq & C\norm{h}_* \left\{ \frac{1}{\gamma_j} + \norm{\mathcal{L}(\tilde{Z}_{0j})}_* \right\} + C\sum_{k=1}^m \abs{d_k}\norm{\mathcal{L}(\tilde{Z}_{0k})}_*\norm{\mathcal{L}(\tilde{Z}_{0j})}_* + \\
& & \quad \quad + C\sum_{l=1}^2\sum_{k=1}^m \abs{e_{lk}}\frac{\norm{\mathcal{L}(\tilde{Z}_{0j})}_{*}}{\gamma_k}
\end{eqnarray*}
where we have used that
$$\displaystyle{\int_{\Omega_s}\frac{\gamma_j}{(\gamma_j^2 +
\abs{y-\xi_j'}^2)^{3/2}} \leq C} \quad \hbox{for all } j .$$
 But
estimate (\ref{App1}) and Claim~\ref{Cl1} imply

\begin{eqnarray}\label{App6}
\abs{d_j}\int_{\Omega_s}\mathcal{L}(\tilde{Z}_{0j})\tilde{Z}_{0j}& \leq & C\frac{\norm{h}_*}{\gamma_j} + C\sum_{k=1}^m \frac{\abs{d_k}\log ^2 s }{ \gamma_j\gamma_k \abs{\log \delta_j}\abs{\log \delta_k}}+
\end{eqnarray}
$$\quad\quad\quad + C\sum_{k\neq j}\abs{d_k}\left|\int_{\Omega_s}\mathcal{L}(\tilde{Z}_{0j})\tilde{Z}_{0k}\right|.$$
We only need to estimate the terms $\displaystyle{\int_{\Omega_s} \mathcal{L}(\tilde{Z}_{0j})\tilde{Z}_{0k}}$, for all $k$. We have the following

\begin{claim}\label{Cl3} If $R$ is sufficiently large,
\begin{equation}\label{App7}
\int_{\Omega_s} \mathcal{L}(\tilde{Z}_{0j})\tilde{Z}_{0j}= \frac{E}{\gamma_j^2\abs{\log\delta_j}}(1+o(1)),
\end{equation}
where $E$ is a possitive constant independent of $s$ and $R$. Besides, if $k\neq j$
\begin{equation}\label{App8}
\int_{\Omega_s} \mathcal{L}(\tilde{Z}_{0j})\tilde{Z}_{0k}= O\left( \frac{\log ^2 s }{\gamma_j\gamma_k \abs{\log\delta_j}\abs{\log\delta_k }} \right).
\end{equation}
\end{claim}

\medskip
Assuming for the moment the validity of  this claim, then replacing (\ref{App7}) and (\ref{App8}) in (\ref{App6}),
we get
\begin{equation}\label{App9}
\frac{\abs{d_j}}{\gamma_j}\leq C\abs{\log\delta_j}\norm{h}_* + C\sum_{k=1}^m \frac{\abs{d_k}}{\gamma_k}\frac{\log^2 s }{\abs{\log\delta_k}},
\end{equation}
and then, $$\abs{d_j}\leq C\gamma_j\abs{\log\delta_j}\norm{h}_*.$$
Finally, using estimate (\ref{App1}), we conclude
$$\abs{e_{ij}}\leq C\gamma_j\log s\norm{h}_*$$ and Claim \ref{Cl2}
holds.
Let us proof Claim \ref{Cl3}. Let us try with the first term (\ref{App7}). We decompose

\begin{eqnarray*}
\int_{\Omega_s}\mathcal{L}(\tilde{Z}_{0j})\tilde{Z}_{0j} & = & O(s\delta_j) +
\int_{\Omega_{2}}\mathcal{L}(\tilde{Z}_{0j})\tilde{Z}_{0j} + \int_{\Omega_{3}}\mathcal{L}(\tilde{Z}_{0j})\tilde{Z}_{0j} + \int_{\Omega_{4}}\mathcal{L}(\tilde{Z}_{0j})\tilde{Z}_{0j}\\
& \equiv & O(s\delta_j) + I_2 + I_3 + I_4.
\end{eqnarray*}
First we estimate $I_3$. From (\ref{5}),

\begin{eqnarray*}
I_3 & = & \int_{\Omega_3}\mathcal{L}(\hat{Z}_{0j})\hat{Z}_{0j}\\
& = & \int_{\Omega_{3,j} \cup \tilde{\Omega}_3}\mathcal{L}(\hat{Z}_{0j})\hat{Z}_{0j} + \sum_{k\neq j}\int_{\Omega_{3,k}}\mathcal{L}(\hat{Z}_{0j})\hat{Z}_{0j}\\
& = & O\left( \frac{1}{R\gamma_j^2 \abs{\log\delta_j}}\right) + O\left( \frac{\log^2 s}{\gamma_j^2 \abs{\log\delta_j}^2}\right).
\end{eqnarray*}
Now we estimate $I_4$. From the estimates in $\Omega_4$,
$\displaystyle{ \abs{I_4}=O\left(
\frac{1}{\gamma_j^2\abs{\log\delta_j}^2}\right)}$.
On the other hand, we have
\begin{equation*}
I_2=\int_{\Omega_{2}}\left\{ \Delta \eta_{1l}(Z_{0j}-\hat{Z}_{0j})
+ 2\nabla\eta_{1j} \cdot \nabla
(Z_{0j}-\hat{Z}_{0j})\right\}\tilde{Z}_{0j} +
O\left(\frac{s\delta}{\gamma_j R^2}\right).
\end{equation*}
Thus integrating by parts the first term above we find
\begin{eqnarray*}
I_2 & = &\int_{\Omega_{2}} \nabla\eta_{1j} \cdot \nabla
(Z_{0j}-\hat{Z}_{0j})\hat{Z}_{0j} -\int_{\Omega_{2}}
\abs{\nabla\eta_{1l}}^2(Z_{0j}-\hat{Z}_{0j})^2  + \\
& & - \int_{\Omega_{2}}
\nabla\eta_{1j} \cdot \nabla\hat{Z}_{0j} (Z_{0j}-\hat{Z}_{0j}) +  \int_{\Omega_{2}}\big\{\eta_{1j}\mathcal{L}(Z_{0j}) + (1-\eta_{1j})\mathcal{L}(\hat{Z}_{0j})\big\} \\
& \equiv & I_{2,a} + I_{2,b} + I_{2,c} + I_{2,d}.
\end{eqnarray*}
Using (\ref{App3}) and (\ref{Lp22a}), we get $\displaystyle{ \abs{\nabla
\hat{Z}_{0l}}}=O\Big(\frac{1}{R^3\gamma_j^2}\Big)$ in $\Omega_2$,
\begin{equation*}
I_{2,b}=O\left(\frac{R}{\gamma_j^2 \abs{\log\delta_j}^2} \right), \quad I_{2,c}=O\left(
\frac{1}{R^2\gamma_j^2 \abs{\log\delta_j}}\right),
\end{equation*}
$$\hbox{ and } \quad I_{2,d}=O\left(\frac{\delta}{R^3\gamma_j^2 \abs{\log\delta_j}}\right).$$
Now, as $\hat{Z}_{0j}=Z_{0j}\left( 1+ O\left(\frac{\gamma_j\delta R}{\abs{\log\delta_j}}\right)\right)$, we conclude
\begin{eqnarray*}
I_{2,a} & = &\frac{1}{\gamma_j^2 \abs{\log\delta_j}}\int_R^{R+1}r\eta_1'(r)\left(\frac{1-r^2}{1+ r^2}\right)(1+o(1))\,dr\\
& = & \frac{E}{\gamma_j^2 \abs{\log\delta_j}}(1+o(1)),
\end{eqnarray*}
where $E$ is a positive constant independent of $s$ and $R$. Thus,
for fixed $R$ large and $s$ small, we obtain (\ref{App7}). The
second result can be established  with similar
 arguments.  \end{proof}

\medskip
Now we can now treat the original linear problem (\ref{Lp5}).

\bigskip
\textbf{Proof of Proposition~\ref{LpTeo1}}. We first establish the
validity of the a priori estimate (\ref{LpTeo1a}) for solutions
$\psi\in L^{\infty}(\Omega)$ of problem (\ref{Lp5}), with $h\in
L^{\infty}(\Omega)$.  Lemma \equ{LpLe4} implies

\begin{equation}\label{Lp40}
\norm{\psi}_{\infty}\leq Cs \Big\{ \norm{h}_{*} + \sum_{i=1}^2
\sum_{j=1}^m \abs{c_{ij}}\norm{\chi_jZ_{ij}}_{*}\Big\},
\end{equation}
but $$\norm{\chi_jZ_{ij}}_{*}\leq C\gamma_j, $$ then, it is
sufficient to estimate the values of the constants $c_{ij}$. To
this end, we multiply the first equation in (\ref{Lp5}) by
$Z_{ij}\eta_{2j}$, with $\eta_{2j}$ the cut-off function
introduced in (\ref{Lp23}), and integrate  by parts to find

\begin{equation}\label{Lp41}
\int_{\Omega_s}\psi\mathcal{L}(Z_{ij}\eta_{2j})=\int_{\Omega_s}hZ_{ij}\eta_{2j}
+ \sum_{k=1}^2 \sum_{l=1}^m c_{kl}\int_{\Omega_s}\eta_{2j}Z_{ij}\chi_l Z_{kl},
\end{equation}
It is easy to see that $\displaystyle{\int_{\Omega_s}h\eta_{2j}Z_{ij}=O(\gamma_j^{-1}\norm{h}_{*})}$. On the other hand we have

\begin{eqnarray*}\label{Lp42}
\mathcal{L}(\eta_{2j}Z_{ij})& = & \Delta\eta_{2j}Z_{ij} + 2\nabla\eta_{2j}\cdot\nabla Z_{ij} + \eta_{2j}\mathcal{L}(Z_{ij})\\
& = & O(\delta^3) + \left\{ W- \frac{8\gamma_j^2}{(\gamma_j^2 + \abs{y-\xi_j'}^2)^2}\right\}\eta_{2j}Z_{ij} \equiv O(\delta^3) + B_j
\end{eqnarray*}
To estimate  $B_j$, we need to split $\supp \eta_{2j}$ into
several pieces. We consider the following subdomains. For a fixed
$j$, we let $$\displaystyle{\hat{\Omega}_{1k}\equiv \left\{
\abs{y-\xi_k'} \leq \frac{1}{2s^{\beta}\delta}\right\}},$$
 for any
$k=1,\dots,m$, and
\begin{equation*}
\hat{\Omega}_{2}\equiv \left\{ \abs{y-\xi_j'} \leq \frac{1}{3\delta},\; \abs{y-\xi_k'} \geq \frac{1}{2s^{\beta}\delta},\; \forall k\right\}.
\end{equation*}
In $\hat{\Omega}_{1j}$, using Lemma \ref{AsmLe1}, $\displaystyle{B_j=O\left( \frac{s\delta\gamma_j}{(\gamma_j^2 + \abs{y-\xi_j'}^2)^{3/2}}\right)}$. In $\hat{\Omega}_{1k}$, $k\neq j$, $$\displaystyle{B_j=O\left( \frac{s^{\beta}\delta\gamma_k^2}{(\gamma_k^2 + \abs{y-\xi_k'}^2)^2}\right)}.$$ Finally, in $\hat{\Omega}_{2}$, $B_j=O(s^K\delta_j^2\delta^3)$, for some constant $K>0$ large. Then,

\begin{equation*}
\left| \int_{\Omega_s} \psi \mathcal{L}(\eta_{2j}Z_{ij})\right|\leq
Cs^{\beta}\delta\norm{\psi}_{\infty}
\end{equation*}
Now, $$\displaystyle{\int_{\Omega_s} \eta_{2j}\chi_j
Z_{ij}Z_{kl}=C\delta_{ik}}$$ and if $j\neq l$, and $s$ is
sufficiently large,
\begin{equation*}
\int_{\Omega_s} \eta_{2j}\chi_l Z_{ij}Z_{kl}=O(\gamma_ls^{\beta}\delta).
\end{equation*}
Using the above estimates in (\ref{Lp41}), we obtain

\begin{equation}\label{Lp43}
\abs{c_{ij}}\leq Cs^{\beta}\delta \norm{\psi}_{\infty} + \frac{C}{\gamma_j} \norm{h}_* + C\sum_{k=1}^2\sum_{l\neq j} \abs{c_{kl}}\gamma_ls^{\beta} \delta
\end{equation}
and then $$\displaystyle{\abs{c_{ij}}\leq Cs^{\beta}\delta
\norm{\psi}_{\infty} + \frac{C}{\gamma_j} \norm{h}_*}.$$ Putting
this estimate in (\ref{Lp40}), we conclude the validity of
(\ref{LpLe4a}).

\medskip




Finally, the a priori estimate implies in particular that the
homogeneous problem has only the trivial solution. A standard
argument involving Fredholm's alternative, see e.g. \cite{KMP},
gives existence. This concludes the proof. \qed

\bigskip
\begin{remark}\label{Re1} {\em
 The operator $T$ is differentiable with  respect to the variables
 $\xi'$. In fact,  computations similar to those used in
 \cite{KMP} yield
 the  estimate
\begin{equation}\label{Re1a}
\norm{\partial_{\xi'}T(h)}_{\infty}\leq Cs^2 \norm{h}_*,\quad
\hbox{ for all }l=1,2;\;k=1,\dots,m.
\end{equation}
Important element in this computation is that $\frac
{1}{\gamma_j}\leq C$, uniformly on $s$.}
\end{remark}

\bigskip

\section{The intermediate nonlinear problem}

In order to solve Problem (\ref{Asm7}) we consider first the
intermediate nonlinear problem.

\begin{eqnarray}\label{NLP1}
\begin{cases} & \mathcal{L}(\psi)= -[R + N(\psi)] + \sum_{i=1}^2 \sum_{j=1}^m c_{ij}\chi_jZ_{ij}, \quad \hbox{ in } \Omega_s, \\
& \psi=0, \quad \hbox{ on } \partial \Omega_s,\\
& \int_{\Omega_s}\chi_jZ_{ij}\psi =0, \hbox{ for all } i=1,2,\;
j=1,\dots , m.
\end{cases}
\end{eqnarray}
For this problem we will prove

\begin{prop}\label{NLPLe1}
Let $\xi\in \mathcal{O}_{s}$. Then, there exists $s_0>0$ and $C>0$ such that for all $s\geq s_0$
the nonlinear problem (\ref{NLP1}) has a unique solution $\psi\in $
which satisfies
\begin{equation}\label{NLP2}
\norm{\psi}_{\infty} \leq C\,s^{2\beta +1} e^{-s/4}.
\end{equation}
Moreover, if we consider the map
$\xi'\in\mathcal{O}_{s}\to\psi\in\mathcal{C}(\overline{\Omega}_s)$,
the derivative $D_{\xi'}\psi$ exists and defines a continuous map
of $\xi'$. Besides
\begin{equation}\label{NLP3}
\norm{D_{\xi'}\psi}_{\infty} \leq C\,s^{2\beta +2} e^{-s/4}.
\end{equation}
\end{prop}

\begin{proof}
In terms of the operator $T$ defined in Proposition \ref{LpTeo1},
Problem (\ref{NLP1}) becomes

\begin{equation*}
\psi = \mathcal{B}(\psi)\equiv -T(N(\psi)+ R).
\end{equation*}
Let us
consider the region
\begin{equation*}
\mathcal{F}\equiv \{ \psi\in
\mathcal{C}(\overline{\Omega}_s) \,\|\, \norm{\psi}_{\infty}\leq s^{2\beta +1}e^{-s/4} \}.
\end{equation*}
From Proposition~\ref{LpTeo1}, $$\norm{\mathcal{B}(\psi)}_{\infty}\leq
C\,s\,\big\{ \norm{N(\psi)}_{*} + \norm{R}_{*}  \big\}, $$ and
Lemma~\ref{AsmLe1} implies $$\norm{R}_{*}\leq Cs^{2\beta+1}e^{-s/4}.$$ Also, from the definition of $N$ in
(\ref{Asm7}), Mean-Value theorem and Lemma~ \ref{AsmLe2} we obtain
$$\norm{N(\psi)}_{*}\leq \norm{W}_{*}\norm{\psi}^2_{\infty}\leq
C\norm{\psi}^2_{\infty}.$$ Hence, if $\psi\in
\mathcal{F}_{\gamma}$, $\norm{\mathcal{B}(\psi)}_{\infty}\leq
Cs^{2\beta + 2} e^{-s/4}$. Along the same way we obtain
$\norm{N(\psi_1)-N(\psi_2)}_{*}\leq C
\max_{i=1,2}\norm{\psi_i}_{\infty}\norm{\psi_1-\psi_2}_{\infty}$,
for any $\psi_1,\psi_2\in \mathcal{F}_{\gamma}$.  Then, we
conclude

\begin{equation*}
\norm{\mathcal{B}(\psi_1)-\mathcal{B}(\psi_2)}_{\infty}\leq
Cs\,\norm{N(\psi_1)-N(\psi_2)}_{*}\leq Cs^{2\beta +2 }e^{-s/4}
\norm{\psi_1-\psi_2}_{\infty}
\end{equation*}
It follows that for all $s$ sufficiently large $\mathcal{B}$ is a
contraction mapping of $\mathcal{F}_{\gamma}$, and therefore a
unique fixed point of $\mathcal{B}$ exists in this region. The
proof of (\ref{NLP3}) is similar to one included in  \cite{KMP}
and we thus omit it.\end{proof}
\medskip


\section{Variational reduction}

We have solved the nonlinear problem (\ref{NLP1}). In order to
find a solution to the original problem (\ref{Asm7}) we need
to find $\mathbf{\xi}$ such that

\begin{equation}\label{VR1}
c_{ij}=c_{ij}(\mathbf{\xi'})=0,\quad \hbox{for all } i,j.
\end{equation}
where $c_{ij}(\mathbf{\xi'})$ are the constants in (\ref{NLP1}).
Problem (\ref{VR1}) is indeed variational: it is equivalent to
finding critical points of a function of $\mathbf{\xi'}$. In fact,
we define the functional for $\mathbf{\xi}\in\mathcal{O}_{s}$:
\begin{equation}\label{VR2}
\mathcal{F}(\mathbf{\xi}) \equiv
J_s[U(\mathbf{\xi})+\hat{\psi}_{\mathbf{\xi}}]
\end{equation}
where $U(\mathbf{\xi})$ is our approximate solution from (\ref{Aa5a}) and $\hat{\psi}_{\mathbf{\xi}}=\psi\big(\frac
{x}{\delta},\frac{\mathbf{\xi}}{\delta}\big)$, $x\in \Omega$, with
$\psi=\psi_{\mathcal{\xi}'}$ the unique solution to problem
(\ref{NLP1}) given by Proposition \ref{NLPLe1}. Then we obtain that
critical points of $\mathcal{F}$ correspond to solutions of
(\ref{VR1}) for large $s$. That is,
\begin{lema}\label{VRLe1} $\mathcal{F}:\mathcal{O}_{s}\to \R$ is of class $\mathcal{C}^1$. Moreover, for all $s$ sufficiently large, if $D_{\mathbf{\xi}}\mathcal{F}(\mathbf{\xi})=0$ then $\mathbf{\xi}$ satisfies (\ref{VR1}).
\end{lema}

\begin{proof} The proof of this fact is standard, see  \cite{KMP}, \cite{EGP} or \cite{Wei}.
Here the estimate found for $D_{\xi'}\psi$ is
used.
\end{proof}

\medskip
The estimates for the solution $\psi_{\xi'}$ for Problem
(\ref{NLP1}) in Proposition \ref{NLPLe1} and  a Taylor expansion
of $\mathcal{F}$ in the expanded domain $\Omega_s$ similar to one
done in \cite{KMP}  give us

\begin{lema}\label{VRLe2} For points $\mathcal{\xi}\in\mathcal{O}_s$ the following expansion holds
\begin{equation}\label{VRLe21}
\mathcal{F}_s(\mathbf{\xi})= J_s[U(\mathbf{\xi})]
+\theta_s(\mathbf{\xi}),
\end{equation}
where $\abs{\theta_s}=O(s^{K}e^{-s/2})$, for some fixed constant $K>0$, uniformly on $s$.
\end{lema}

\medskip


\section{The Proof of Theorem \ref{teo1}}

\medskip We consider the set
\begin{equation}S =\{ x\in \Lambda \,| \, \phi_1(x)=1 \}.
\label{S}\end{equation}
The result Theorem~\ref{teo1}  is a direct
consequence of the following more precise result.
\begin{teo}\label{T1}
Given any positive integer $m$ there exists $s_0>0$ sufficiently
large such that problem $(\ref{In3})$ has a solution $u_s$
positive in $\Omega$ of the form

\begin{equation}\label{T1a}
u_s(x)=U(\xi^s)+ \tilde{\psi}_s,
\end{equation}
which possesses exactly $m$ local maximum points $\xi_1^s,\dots ,
\xi_m^s\in\Lambda$, satisfying that as $s\to \infty$

{\rm (i) } $ \dist(\xi_j^s , S) \to 0 $  and
$\abs{\xi_i^s-\xi_j^s}\geq \frac{1}{s^{m(m+1)}}$ if $i\neq j$;

{\rm (ii) }  $\norm{\tilde{\psi}_s}_{\infty} \to 0 $.
\end{teo}

\medskip
The construction actually yields $1-\phi_1(\xi_j^s) > s^{-\frac
12}$. Thus if $S$ is just constituted by a non-degenerate maximum
point $\bar x$ we will have $|\xi_j^s - \bar x| \le C s^{-\frac 14
}$.

\begin{proof}
 According to
Lemma~\ref{VRLe1}, $U(\xi^s) + \hat{\psi}_{\xi^s}$ is a solution
of problem (\ref{In3}) if  $\xi^s\in \mathcal{O}_s$ is a critical
point of the functional $\mathcal{F}$ defined in \equ{VR2}.
 We recall in particular that $\norm{\tilde \psi_{\xi^s } }_{\infty} \to 0$ as predicted by
estimate (\ref{NLP2}). It thus suffices to establish that
$\mathcal{F}$ attaints its maximum value in $\mathcal{O}_s$ for
all sufficiently large $s$, for which we will see
\begin{equation}
\sup_{\xi \in \partial \mathcal{O}_s} \mathcal{F}(\xi)\, < \,
\sup_{\xi\in \mathcal{O}_s} \mathcal{F}(\xi) .
\label{max}\end{equation} First we obtain a lower bound  for
$\sup_{\xi\in \mathcal{O}_s} \mathcal{F}(\xi)$
 Let us fix a point $\bar{x}\in S$ and
set
\begin{equation*}\label{PT1a}
\xi ^0_j\equiv \bar{x}+ \frac {1}{\sqrt{s}}\hat{\xi}_j,
\end{equation*}
where $\hat{\xi}=(\hat{\xi}_1,\dots,\hat{\xi}_m )$ is a $m$-regular
polygon in $\R^2$. Clearly $\xi^0\in \mathcal{O}_{s}$ because $\phi_1(\xi_j^0)=1+O(s^{-1})$. Then

\begin{eqnarray*}\label{PT2}
\sup_{\xi\in \mathcal{O}_s} \mathcal{F}(\xi) & \ge &  J_s(U(\xi^0)) + \theta_s(\xi^0)\\
&  = & 8\pi\Big\{ \sum_{i\neq j}^m 2\log\abs{\xi_i^0 -\xi_j^0}
+ s \sum_{j=1}^m \phi_1(\xi_j^0) \Big\} + O(1)\\
& \geq & 8\pi m \Big\{ -(m-1)\log s + s \Big\} + O(1)
\end{eqnarray*}
Then,
\begin{equation}\label{B1}
\sup_{\xi\in \mathcal{O}_s} \mathcal{F}(\xi)  \geq 8\pi m s -8\pi
m(m-1)\log s + O(1).
\end{equation}
Next we estimate from above
 $\mathcal{F}(\xi)$ for
 $\xi\in\partial\mathcal{O}_{s}$. Then, there are two
possibilities: either (1) there exist indices $i_0,j_0$, $i_0\neq
j_0$ such that $\abs{\xi_{i_0}-\xi_{j_0}}=s^{-\beta}$, or (2)
there exists $i_0$ such that
$1-\phi_1(\xi_{i_0})=\frac{1}{\sqrt{s}}>0$.

\medskip
In the first case, we have the following upper bound
\begin{eqnarray}\label{PT3}
\mathcal{F}(\xi)
 \leq  8\pi\big\{ -2\beta \log s+ s \sum_{j=1}^m \phi_1(\xi_j) \big\} + O(1)
 \leq  8\pi m \big\{ s  -\frac 2m \beta\log s\big\} + O(1).
\end{eqnarray}

\medskip
\noindent
In the second case, $1-\phi_1(\xi_{i_0}^s)\leq
\frac{1}{2\sqrt{s}}$.  Then

\begin{eqnarray}\label{PT5}
 \mathcal{F} (\xi)  \leq  8\pi\Big\{ O(\log s)+ s \Big(1-\frac{1}{2\sqrt{s}} + (m-1)\Big)\Big\} 
 \leq  8\pi s\big( m- \frac{1}{2s^{1/2}} \big) + O(\log s) \, .
\end{eqnarray}
 At this
point we make the election $\beta> m^2 +m$ in the definition of
${\mathcal O}_s$. Relation \equ{max} immediately follows from
combining estimates \equ{PT2}, \equ{PT3}, \equ{PT5} and taking $s$
sufficiently large. This finishes the proof.\end{proof}

\bigskip

 \centerline{\bf Acknowledgement}

\medskip
This work has been partly supported by grants Fondecyt 1030840 and
FONDAP, Chile.

\end{document}